\documentclass[11pt]{amsart}

\usepackage{amsmath}
\usepackage{latexsym}
\usepackage{amssymb}
\usepackage{amsfonts}
\usepackage{graphics}
\usepackage{CJK}
\usepackage{amscd}
\usepackage{arydshln}
\usepackage{latexsym,amsmath,amsfonts,amssymb,stmaryrd,bbding,bbm}
\usepackage{fancyhdr,color,graphicx}
\usepackage{mathrsfs}
\usepackage[centertags]{amsmath}
\usepackage{amssymb}
\usepackage{amsthm}
\usepackage{newlfont}

\include{PDF}

\newtheorem{theorem}{Theorem}[section]%
\newtheorem{lemma}[theorem]{Lemma}%
\newtheorem{corollary}[theorem]{Corollary}
\newtheorem{{}}[theorem]{}%

\makeatletter\def\equfill@{arrowfill@Relbar\Relbar\Relbar}

\newcommand{\equfill}[2]{\ext@arrow 0395\equfill@{#1}{#2}}
\makeatother
\def\f{\noindent}

\def\qed{\hfill $\Box$}

\def\demo{{\bf Proof}\hskip10pt}

\begin{document}

\title{ \bf  $p$-Sylowizers and $p$-nilpotency of finite groups}
\thanks{This work was supported by the National Natural Science Foundation of China(Grant N.12261022, 11871360 and 12201495)}
\author{ Yaxin Gao$^1$, Xianhua Li$^2$, Donglin Lei$^3$}

\address{$^1$Department of Mathematics, Changzhou University, Changzhou 213164, P. R. China}
\address{$^2$School of Mathematical Science, Guizhou Normal University, Guiyang 550025, P. R. China}
\address{$^3$Department of Mathematics, Northwest University, Xi'an 710127, P. R. China}
\email{yxgaomath@163.com}
\maketitle

\begin{abstract} In this paper, we investigate the structure of finite group $G$ by assuming that the intersections between $p$-sylowizers of some $p$-subgroups of $G$ and $O^p(G)$ are $S$-permutable in $G$. We obtain some criterions for $p$-nilpotency of a finite group.

\medskip
\vspace{0.2cm}
 \f {\bf Key words:} finite groups, $p$-sylowizers, subgroups, permutable, Sylow subgroups
\vspace{0.2cm}

\f {\bf AMS Subject Classifications:} 20D10, 20D20
\end{abstract}

\section{Introduction}

All groups in this paper are finite, $G$ is always a finite group. $\pi$ denotes a set of primes, $G_{\pi}$ means a Hall $\pi$-subgroup of $G$, $O_{\pi}(G)$ is the largest normal $\pi$-subgroup of $G$, and $O^{\pi}(G)$ is the subgroup generated by all $\pi'$-elements of $G$. We use conventional notions and notations, as in \cite{HB} and \cite{KT}.

The concept of $p$-sylowizers was introduced by Gasch$\ddot{u}$tz  \cite{G}: A subgroup
$S$ of a group $G$ is called a $p$-sylowizer of a $p$-subgroup $R$ in $G$ if $S$ is maximal in $G$ with respect to having $R$ as its Sylow $p$-subgroup. Recall that
a subgroup $H$ of a group $G$ is said to be $S$-permutable (or $S$-quasinormal, $\pi$-quasinormal in some literatures) in $G$ if $HG_{p}=G_{p}H$ for all Sylow subgroups $G_{p}$ of $G$ (see \cite{O}). In \cite{LL}, Lei and Li obtained some new characterizations of $p$-nilpotent and supersolvable groups by the permutability of the $p$-sylowizers of some $p$-subgroups.

Recently, Li and Zhang \cite{LZ} investigated the $p$-supersolvability of a finite group by analyzing the intersections between $p$-sylowizers of some $p$-subgroups with given order and $O^p(G)$. One of their main theorems is as follows:

\begin{theorem}  (\cite[Theorem 1.1]{LZ}) Let $p$ be a prime dividing the order of $G$, $P\in Syl_{p}(G)$, and $d$ be a divisor of $|P|$ with $1\leq d<|P|$. Suppose that for every subgroup $H$ of $P$ with order $d$, $H$ has a $p$-sylowizer $S$ in $G$ with $S\cap O^p(G)\unlhd O^p(G)$. Then $G$ is $p$-supersolvable.
\end{theorem}

Continuing the work of Li and Zhang, Yu, Du and Xu \cite{Y2} proved the following theorem.

\begin{theorem}  (\cite[Theorem 1.3]{Y2}) Let $p$ be a prime dividing the order of $G$, $P\in Syl_{p}(G)$, and $d$ be
a divisor of $|P|$ with $1\leq d<|P|$. Suppose that for every normal subgroup $H$
of $P$ with order $d$, $H$ has a $p$-sylowizer $S$ in $G$ with $S\cap O^p(G)\unlhd O^p(G)$. Then $G$ is $p$-nilpotent.
\end{theorem}


Continuing the study above, we investigate the influence of the intersections between $p$-sylowizers of every normal subgroup $H$ of a Sylow $p$-subgroup of $G$ with order $d$ and $O^p(G)$ on the structure of a finite group $G$ by assuming that the intersections are $S$-permutable in $G$. We obtain some new characterizations of $p$-nilpotent groups (see Theorem 3.1, 3.4). In this paper, we also give a characterization on a $p$-nilpotent group $G$ by the properties of a subgroup $P$ satisfing $G_{p}'\leq P \leq \Phi(G_{p})$ for $G_{p}\in Syl_{p}(G)$ (see Theorem 3.6).

As a generalization of $S$-permutable subgroup, the concept of $\mathfrak{Z}$-permutable was introduced. Let $G$ be a finite group and $\mathfrak{Z}$ a complete set of Sylow subgroups of $G$, that is, $\mathfrak{Z}$ contains exactly one Sylow $p$-subgroup of $G$ for each prime $p$ dividing the order of $G$. Then a subgroup $H$ of $G$ is said to be $\mathfrak{Z}$-permutable in $G$ if $H$ permute with all subgroups in $\mathfrak{Z}$ of $G$ (see \cite{LLW}). In this paper, we also give a characterization of a $p$-nilpotent group by assuming that $G_{p}\in \mathfrak{Z}$ and $G_{p}$ has a chain of subgroups $1=P_{0}< P_{1}< \cdots < P_{n-1}< P_{n}=G_{p}$ such that $|P_i:P_{i-1}|=p$ and every $p$-sylowizer of $P_{i}$ in $G$ is $\mathfrak{Z}$-permutable in $G$ for $i=1, \cdots, n$ (see Theorem 3.7).

\section{Preliminaries}
In this section, for the sake of convenience, we list some known results which will be used in the proofs of section $3$ in this paper.

\begin{lemma}\label{L1}(\cite[Lemma 2.3]{LL})
Let $R$ be a $p$-subgroup of $G$ and $S$ a $p$-sylowizer of $R$ in $G$. If $S$ is $S$-permutable in $G$, then $O^p(G)\leq S$. In particular, $S=RO^p(G)$ is the unique $p$-sylowizer of $R$ in $G$.
\end{lemma}

\begin{lemma}\label{L2}(\cite[Lemma 2.2]{LL})
Let $R$ be a $p$-subgroup of $G$. Assume that $N$ is a normal subgroup of $G$ and $R$ is a Sylow $p$-subgroup of $RN$. Then $S$ is a $p$-sylowizer of $R$ in $G$ if and only if $S/N$ is a $p$-sylowizer of $RN/N$ in $G/N$.
\end{lemma}

\begin{lemma}\label{L3}(\cite{O})
Suppose that $H$ is an $S$-permutable subgroup of $G$ and $N\trianglelefteq G$. Then the following holds:

(1) If $K \leq G$, then $H \cap K$ is $S$-permutable in $K$.

(2) $HN$ and $H \cap N$ are $S$-permutable in $G$, $HN/N$ is $S$-permutable in $G/N$.

(3) $H$ is subnormal in $G$.

\end{lemma}

\begin{lemma}\label{L4}(\cite[Lemma 2.1]{MA2})
Let $G$ be a group and $\mathfrak{Z}$ be a complete set of Sylow subgroups of $G$ and $N$ be a normal subgroup of $G$.

(1) $\mathfrak{Z}\cap N$ and $\mathfrak{Z}N/N$ are complete sets of Sylow subgroups of $N$ and $G/N$, respectively.

(2) If $U$ is a $\mathfrak{Z}$-permutable subgroup of $G$, then $UN/N$ is $\mathfrak{Z}N/N$-permutable. If $U$ is contained in $N$, then $U$ is $\mathfrak{Z}\cap N$-permutable.
\end{lemma}

\begin{lemma}\label{L5} Let $H$ be a $p$-group of a group $G$ and $N$ a normal subgroup of $G$. Assume that every $p$-sylowizer of $H$ in $G$ is $\mathfrak{Z}$-permutable in $G$. If $N$ is a $p'$-group or $N\leq H$, then every $p$-sylowizer of $HN/N$ in $G/N$ is $\mathfrak{Z}N/N$-permutable $G/N$.
\end{lemma}

\demo Assume that $S/N$ is a $p$-sylowizer of $HN/N$ in $G/N$. Then by \cite[Hilfssatz]{G}, $S$ is a $p$-sylowizer of $H$ in $G$. Since $S$ is $\mathfrak{Z}$-permutable in $G$, $S/N$ is $\mathfrak{Z}N/N$-permutable in $G/N$ by Lemma \ref{L4}. \qed

\begin{lemma}\label{L6}(\cite[Lemma 2.1]{LL})
Let $H$ be a $p$-subgroup of $G$ and $K$ a subgroup of $G$ satisfying $H \leq K \leq G$. Assume that $T$ is a $p$-sylowizer of $H$ in $K$. Then there is a $p$-sylowizer $S$ of $H$ in $G$ such that $T=S\cap K$.
\end{lemma}

\begin{lemma}\label{L7}(\cite[Lemma A]{P})
If $H$ is $S$-permutable in a group $G$ and $H$ is a $p$-group for some prime $p$, then $O^p(G)\leq N_{G}(H)$.
\end{lemma}

\begin{lemma}\label{L8}(\cite[Lemma 2.4]{WW2})
Let $H$ be a normal subgroup of a group $G$ such that $G/H$ is $p$-nilpotent and let $P$ be a Sylow $p$-subgroup of $H$, where $p$ is a prime divisor of $|G|$. If $|P|\leq p^2$ and one of the following conditions holds, then $G$ is $p$-nilpotent:

(1) $(|G|, p-1)=1$ and $|P|\leq p$;

(2) $G$ is $A_{4}$-free if $p=min$ $\pi(G)$;

(3) $(|G|, p^2-1)=1$.
\end{lemma}

\section{Main Results}

\begin{theorem} \label{Th2} Let $G$ be a group and $G_{p}\in Syl_{p}(G)$ and $d$ be a divisor of $|G_{p}|$ such that $1\leq d< |G_{p}|$, where $p=min$ $\pi(G)$. Then $G$ is $p$-nilpotent if and only if for every normal subgroup $H$ of $G_{p}$ with $|H|=d$ and $4$ (if $d=p=2$ and $G_{p}$ is non-abelian), every $p$-sylowizer $S$ of $H$ in $G$ satisfies that $S\cap O^p(G)$ is $S$-permutable in $G$.
\end{theorem}

\demo  Assume that $G$ is $p$-nilpotent. Then there is a normal $p$-complement $K$ of $G$ such that $G=G_{p}\ltimes K$. Clearly, $O^p(G)=K$, and for any normal subgroup $H$ of $G_{p}$ with $|H|=d$, $HK$ is the unique $p$-sylowizer of $H$ in $G$ by Lemma \ref{L1}. Thus we only need to prove the sufficiency. Let $G$ be a counterexample with minimal order.

(1) $O_{p'}(G)=1$.

Assume that $O_{p'}(G)\neq 1$. Let $K/O_{p'}(G)$ be any normal subgroup of $G_{p}O_{p'}(G)/O_{p'}(G)$ with $|K/O_{p'}(G)|=d$, then $K=HO_{p'}(G)$, where $H$ is a normal subgroup of $G_{p}$ of order $d$. Let $S/O_{p'}(G)$ be any $p$-sylowizer of $HO_{p'}(G)/O_{p'}(G)$ in $G/O_{p'}(G)$. Then by Lemma \ref{L2},  $S$ is $p$-sylowizer of $H$ in $G$. Then by the hypothesis and Lemma \ref{L3}(3), $S/O_{p'}(G)\cap O^p(G/O_{p'}(G))=(S\cap O^p(G))/O_{p'}(G)$ is $S$-permutable in $G/O_{p'}(G)$. This shows that $G/O_{p'}(G)$ satisfies the hypothesis of the theorem. Hence $G/O_{p'}(G)$ is $p$-nilpotent by the minimal choice of $G$, which implies that $G$ is $p$-nilpotent, a contradiction.

(2) $|O^p(G)_{p}|\leq d$, where $O^p(G)_{p}\in Syl_{p}(O^p(G))$.

Assume that $O^p(G)=G$. Then every $p$-sylowizer $S$ of every normal subgroup $H$ of $G_{p}$ with $|H|=d$ satisfies that $S=S\cap O^p(G)$ is $S$-permutable in $G$. By Lemma \ref{L1}, we have $S=HO^p(G)=G$. It implies that $|G_{p}|=|H|=d$, a contradiction. Hence $O^p(G)<G$. Now, assume that $|O^p(G)_{p}|>d$. Then we can pick a normal subgroup $H$ of $G_{p}$ such that $H\leq O^p(G)$ and $|H|=d$. Let $S$ be any $p$-sylowizer of $H$ in $G$. Then $S\leq O^p(G)$. By the hypothesis and Lemma \ref{L1}, $O^p(G)\leq S$. It follows that $S=O^p(G)$, which contradicts $|O^p(G)_{p}|>d$. Thus (2) holds.

(3) $G/N$ is $p$-nilpotent, where $N$ is the unique minimal normal subgroup of $G$ contained in $O^p(G)$.

Let $N$ be a minimal normal subgroup of $G$ contained in $O^p(G)$. Then $|N_{p}|\leq d$ by (2), where $N_{p}\in Syl_{p}(N)$. Let $K/N$ be any normal subgroup of $G_{p}N/N$ with $|K/N|=\frac{d}{|N_{p}|}$. Then $K/N=HN/N$, where $H$ is a normal subgroup of $G_{p}$ with $|H|=d$ and $N_{p}\leq H$. Let $S/N$ be any $p$-sylowizer of $K/N$ in $G/N$. By Lemma \ref{L2}, $S$ is $p$-sylowizer of $H$ in $G$.  Then by Lemma \ref{L3}(2), every $p$-sylowizer $S/N$ of $K/N$ in $G/N$ satisfies that $S/N\cap O^p(G/N)$ is $S$-permutable in $G/N$. This shows that $G/N$ satisfies the hypothesis of theorem and so $G/N$ is $p$-nilpotent by the minimal choice of $G$. Since the formation of all $p$-nilpotent groups is saturated, $N$ is the unique minimal normal
subgroup of $G$ contained in $O^p(G)$ and $N\nleq \Phi(G)$.

(4) $G$ is solvable.

Assume that $G$ is non-solvable. Then $p=2$. If $G_{p}$ is cyclic, certainly $G$ is $p$-nilpotent, a contradiction. Thus we may assume that $G_{p}$ is non-cyclic. It is easy to see that there exists a maximal subgroup $P$ of $G_{p}$ such that $N_{p}\nleq P$ by Tate's theorem (see \cite[Satz IV. 4.7]{HB}), where $N_{p}\in Syl_{p}(N)$. Since $P\unlhd G_{p}$, we can choose a subgroup $H$ of $P$ with $|H|=d$ such that $H\unlhd G_{p}$ and $N_{p}\nleq H$. Clearly, $|HN|_{p}>d$. Assume that $HN<G$. Obviously, $HN_{p}$ is a Sylow $p$-subgroup of $HN$. Now, let $H_{i}$ be any normal subgroup of $HN_{p}$ with $|H_{i}|=d$ and $T_{i}$ be any $p$-sylowizer of $H_{i}$ in $HN$. By Lemma \ref{L6}, there is a $p$-sylowizer $S_{i}$ of $H_{i}$ in $G$ such that $T_{i}=S_{i} \cap HN$. Since $T_{i}\cap O^p(HN)=S_{i}\cap HN\cap O^p(HN)\cap O^p(G)$, then by the hypothesis and Lemma \ref{L3}(1) and (2), $T_{i}\cap O^p(HN)$ is $S$-permutable in $HN$. This shows that $HN$ satisfies the hypothesis and so $HN$ is $p$-nilpotent by the minimal choice of $G$. It implies that $N$ is solvable. Since $G/N$ is $p$-nilpotent by (3), $G$ is solvable, a contradiction. Hence $HN=G$, then $SN=G$ for any $p$-sylowizer $S$ of $H$ in $G$. Since $H\unlhd G_{p}$, there is a maximal subgroup $H_{1}$ of $H$ such that $H_{1}\unlhd G_{p}$. We choose a subgroup $N_{1}$ of $N_{p}$ such that $N_{1}\unlhd G_{p}$ and $|N_{1}|=p$ and construct the subgroup $L=H_{1}\times N_{1}$. Then $L\unlhd G_{p}$ and $|L|=d$. Clearly, $N_{p}\nleq L$ and $S_{0}\cap O^p(G)\neq 1$ for some $p$-sylowizer $S_{0}$ of $L$ in $G$.
By the hypothesis and Lemma \ref{L3}(3), $S_{0}\cap O^p(G)$ is subnormal in $G$. Then by \cite[A, Theorem 14.3]{KT}, we have $N\leq N_{G}(S_{0}\cap O^p(G))$ and so $G=S_{0}N\leq N_{G}(S_{0}\cap O^p(G))$. Thus $S_{0}\cap O^p(G)\unlhd G$. By the minimality and uniqueness of $N$, we have $N\leq S_{0}\cap O^p(G)$. It follows that $N_{p}\leq L$, a contradiction.

(5) Final contradiction.

By (1) and (4), we get that $N$ is a $p$-group. Since $G$ is not $p$-nilpotent and $G/N$ is $p$-nilpotent by (3), then by Lemma \ref{L8}, we have $|N|>p$. So $p<|N|\leq d<|G_{p}|$. Since $N\nleq \Phi(G)$ by (3), $N\nleq \Phi(G_{p})$. Hence $G_{p}$ has a maximal subgroup $P$ such that $N\nleq P$ and so $1<|P\cap N|<|N|\leq d\leq |P|$. Thus there is a normal subgroup $P_{1}$ of $G_{p}$ of order $d$ such that $P\cap N<P_{1}\leq P$. Let $S$ be any $p$-sylowizer of $P_{1}$ in $G$. Then by the hypothesis and Lemma \ref{L3}(2), $P_{1}\cap N=S\cap N=S\cap O^p(G)\cap N$ is $S$-permutable in $G$. By Lemma \ref{L7}, $P_{1}\cap N\unlhd O^p(G)$. Since $P_{1}\cap N\unlhd G_{p}$, $P_{1}\cap N\unlhd G$. It is clear that $P_{1}\cap N>1$. Then by the minimality of $N$, we have $P_{1}\cap N=N$ and so $N\leq P_{1}\leq P$, a contradiction. This contradiction completes the proof. \qed

\vspace{5pt}

Let $G_{p}$ be a Sylow $p$-subgroup of a group $G$, $H$ be a normal subgroup of $G_{p}$ and $S$ be a $p$-sylowizer of $H$ in $G$. It is not difficult to prove that the condition `` $S\cap O^p(G)\unlhd O^p(G)$ '' is equivalent to `` $S\cap O^p(G)\unlhd G$ ''. Since $S\cap O^p(G)\unlhd G$ implies that $S\cap O^p(G)$ is $S$-permutable in $G$, by Theorem \ref{Th2}, we have the following corollary.

\begin{corollary} \label{Co3}
Let $G$ be a group and $G_{p}\in Syl_{p}(G)$ and $d$ be a divisor of $|G_{p}|$ such that $1\leq d< |G_{p}|$, where $p=min$ $\pi(G)$. Then $G$ is $p$-nilpotent if and only if for every normal subgroup $H$ of $G_{p}$ with $|H|=d$ and $4$ (if $d=p=2$ and $G_{p}$ is non-abelian), every $p$-sylowizer $S$ of $H$ in $G$ satisfies that $S\cap O^p(G)\unlhd O^p(G)$.
\end{corollary}

\begin{corollary} \label{Co4}
Let $G$ be a group and $p=min$ $\pi(G)$. Then $G$ is $p$-nilpotent if and only if $G$ has a normal subgroup $N$ such that $G/N$ is $p$-nilpotent and for every normal subgroup $H$ of $N_{p}$ with $|H|=d$ and $4$ (if $d=p=2$ and $N_{p}$ is non-abelian), every $p$-sylowizer $S$ of $H$ in $G$ satisfies that $S\cap O^p(G)$ is $S$-permutable in $G$, where $N_{p}\in Syl_{p}(N)$ and $d$ is a divisor of $|N_{p}|$ such that $1\leq d<|N_{p}|$.
\end{corollary}

\demo  The necessity is evident. We only need to prove the sufficiency. By the hypothesis, $G/N$ is $p$-nilpotent. Let $T/N$ be a normal $p$-complement of $G/N$. Then $N_{p}$ is a Sylow $p$-subgroup of $T$. Let $H$ be any normal subgroup of $N_{p}$ with $|H|=d$ and $K$ be any $p$-sylowizer of $H$ in $T$. Then by Lemma \ref{L6}, $K=S\cap T$ for some $p$-sylowizer $S$ of $H$ in $G$. By the hypothesis and Lemma \ref{L3}, $K\cap O^p(T)=S\cap T\cap O^p(T)\cap O^p(G)$ is $S$-permutable in $T$. This shows that $T$ satisfies the hypothesis of Theorem \ref{Th2} and so $T$ is $p$-nilpotent. Let $T_{p'}$ be a normal $p$-complement of $T$, then $T_{p'}\unlhd G$ and so $G$ is $p$-nilpotent.  \qed

\begin{theorem} \label{Th5}  Let $G$ be a group and $G_{p}\in Syl_{p}(G)$ and $d$ be a divisor of $|G_{p}|$ such that $1\leq d< |G_{p}|$, where $p$ is an odd prime and a divisor of $|G|$. Then $G$ is $p$-nilpotent if and only if $N_{G}(G_{p})$ is $p$-nilpotent and for every normal subgroup $H$ of $G_{p}$ with $|H|=d$, every $p$-sylowizer $S$ of $H$ in $G$ satisfies that $S\cap O^p(G)$ is $S$-permutable in $G$.
\end{theorem}

\demo  The necessity is evident, we only need to prove the sufficiency. Let $G$ be a counterexample with minimal order. Similar to the proof of Theorem \ref{Th2}, it is easy to get (1) and (2).

(1) $O_{p'}(G)=1$.

(2) $|O^p(G)_{p}|\leq d$, where $O^p(G)_{p}\in Syl_{p}(O^p(G))$.

(3) If $G_{p}\leq L<G$, then $L$ is $p$-nilpotent.

Clearly, $N_{L}(G_p)\leq N_{G}(G_{p})$ is $p$-nilpotent. Let $H$ be any normal subgroup of $G_{p}$ with $|H|=d$ and $T$ be any $p$-sylowizer of $H$ in $L$. Then by Lemma \ref{L6}, there is a $p$-sylowizer $S$ of $H$ in $G$ such that $T=S\cap L$. By the hypothesis and Lemma \ref{L3}, $T\cap O^p(L)$ is $S$-permutable in $L$. This shows that $L$ satisfies the hypothesis of the theorem. Hence $L$ is $p$-nilpotent by the minimality of $G$.

(4) Let $O^p(G)_{p}=O^p(G)\cap G_{p}$ be a Sylow $p$-subgroup of $O^p(G)$. Then $O^p(G)_{p}$ has a maximal subgroup $P_{2}$ such that $P_{2}$ is $S$-permutable in $G$.

It is easy to see that $O^p(G)$ is not $p$-nilpotent. Then by Tate's theorem (see \cite[Satz IV. 4.7]{HB}), we have $O^p(G)_{p}=O^p(G)\cap G_{p}\nleq \Phi(G_{p})$. It implies that there is a maximal subgroup $P$ of $G_{p}$ such that $O^p(G)_{p}\nleq P$. By (2), $1<|O^p(G)_{p}|\leq d<|G_{p}|$. Hence $|P\cap O^p(G)_{p}|<|O^p(G)_{p}|\leq d\leq|P|$. Then $G_{p}$ has a normal subgroup $P_{1}$ of order $d$ such that $P\cap O^p(G)_{p}<P_{1}\leq P$. Thus $P_{1}\cap O^p(G)_{p}=P\cap O^p(G)_{p}$ and $|O^p(G)_{p}: P_{1}\cap O^p(G)_{p}|=p$. This shows that $O^p(G)_{p}\nleq P_{1}$. Let $S$ be any $p$-sylowizer of $P_{1}$ in $G$. By the hypothesis, $S\cap O^p(G)$ is $S$-permutable in $G$, hence $G_{p}(S_{0}\cap O^p(G))\leq G$. If $G_{p}(S_{0}\cap O^p(G))=G$ for some $p$-sylowizer $S_{0}$, then $O^p(G)=O^p(S_{0}\cap O^p(G))\leq S_{0}\cap O^p(G)\leq S_{0}$. It follows that $O^p(G)_{p}\leq P_{1}$, a contradiction. This contradiction shows that $G_{p}(S\cap O^p(G))<G$ for any $p$-sylowizer $S$. Then by (3), $G_{p}(S\cap O^p(G))$ is $p$-nilpotent, so does $S\cap O^p(G)$. Since $O_{p'}(G)=1$ and $S\cap O^p(G)\lhd\lhd G$ by Lemma \ref{L3}(3), $P_{1}\cap O^p(G)=S\cap O^p(G)$.  Let $P_{2}=P_{1}\cap O^p(G)_{p}$. Hence $|O^p(G)_{p}:P_{2}|=p$ and $P_{2}=P_{1}\cap O^p(G)$ is $S$-permutable in $G$. Moreover, $P_{2}\leq O_{p}(G)$.

(5) Final contradiction.

Let $O^p(G)_{p}=O^p(G)\cap G_{p}$ be a Sylow $p$-subgroup of $O^p(G)$. Assume that $N_{G}(O^p(G)_{p})<G$. By (4), $O^p(G)_{p}$ has a maximal subgroup $P_{2}$ such that $P_{2}$ is $S$-permutable in $G$ and $P_{2}\leq O_{p}(G)$. Then by Lemma \ref{L7}, $P_{2}\unlhd O^p(G)$. It implies that $P_{2}=O_{p}(O^p(G))$. Clearly, $G_{p}\leq N_{G}(O^p(G)_{p})$. Then by (3), $N_{G}(O^p(G)_{p})$ is $p$-nilpotent and so $N_{O^p(G)}(O^p(G)_{p})$ is $p$-nilpotent. We write $A=O^p(G)_{p}/P_{2}$ and $B=O^p(G)/P_{2}$. It is clear that $N_{B}(A)\cong N_{O^p(G)}(O^p(G)_{p})/P_{2}$ is $p$-nilpotent and $A\in Syl_{p}(B)$ and $|A|=p$. By Burnside's theorem, $B$ is $p$-nilpotent. It follows that $O^p(B)<B$, so $O^p(O^p(G))<O^p(G)$, a contradiction. Now assume that $N_{G}(O^p(G)_{p})=G$, that is $O^p(G)_{p}\unlhd G$. Then by Schur-Zassenhaus's theorem, $O^p(G)$ has a Hall $p'$-subgroup $K$. By Frattini's argument, we have $G=O^p(G)N_{G}(K)$ and so $G_{p}=O^p(G)_{p}N_{G_{p}}(K)$. Clearly, $O^p(G)$ is non-$p$-nilpotent, then $O^p(G)_{p}\nleq N_{G_{p}}(K)$ and so $N_{G_{p}}(K)<G_{p}$. It implies that there is a maximal subgroup $P_{0}$ of $G_{p}$ such that $N_{G_{p}}(K)\leq P_{0}$. Then $G_{p}=O^p(G)_{p}P_{0}$ and $|O^p(G)_{p}:O^p(G)_{p}\cap P_{0}|=p$, and so $|O^p(G)_{p}\cap P_{0}|<|O^p(G)_{p}|\leq d\leq |P_{0}|$. Hence $G_{p}$ has a normal subgroup $P$ of order $d$ such that $O^p(G)_{p}\cap P_{0}<P$. Then $O^p(G)_{p}\cap P=O^p(G)_{p}\cap P_{0}$ and so $|O^p(G)_{p}:O^p(G)_{p}\cap P|=p$. Moreover, $O^p(G)_{p}\nleq P$. By a similar argument as (4), we have that $P\cap O^p(G)=S\cap O^p(G)$ is $S$-permutable for any $p$-sylowizer $S$ of $P$ in $G$. By Lemma \ref{L7}, $O^p(G)\leq N_{G}(P\cap O^p(G))$. Since $P\cap O^p(G)=P\cap G_{p}\cap O^p(G)=P\cap O^p(G)_{p}\unlhd G_{p}$, $P\cap O^p(G)\unlhd G$. Clearly, $1<|P\cap O^p(G)|<d$. Consider the group $G/(P\cap O^p(G))$. By Lemma \ref{L2} and Lemma \ref{L3}, it is easy to see that every $p$-sylowizer $S/(P\cap O^p(G))$ of every normal subgroup $H/(P\cap O^p(G))$ of $G_{p}/(P\cap O^p(G))$ with $\frac{d}{|P\cap O^p(G)|}$ satisfies that $S/(P\cap O^p(G))\cap O^p(G/(P\cap O^p(G)))=(S\cap O^p(G))/(P\cap O^p(G))$ is $S$-permutable in $G/(P\cap O^p(G))$ and $N_{G/(P\cap O^p(G))}(G_{p}/(P\cap O^p(G)))\cong N_{G}(G_{p})/(P\cap O^p(G))$ is $p$-nilpotent. This shows that $G/(P\cap O^p(G))$ satisfies the hypothesis and so $G/(P\cap O^p(G))$ is $p$-nilpotent by the minimality of $G$. It implies that $O^p(G)/(P\cap O^p(G))$ is $p$-nilpotent. Hence $O^p(O^p(G)/(P\cap O^p(G)))<O^p(G)/(P\cap O^p(G))$ and so $O^p(O^p(G))<O^p(G)$, a contradiction. This contradiction completes the proof. \qed

\begin{corollary}  Let $G$ be a group and $G_{p}\in Syl_{p}(G)$ and $d$ be a divisor of $|G_{p}|$ such that $1\leq d< |G_{p}|$, where $p$ is an odd prime and a divisor of $|G|$. Then $G$ is $p$-nilpotent if and only if $N_{G}(G_{p})$ is $p$-nilpotent and for every normal subgroup $H$ of $G_{p}$ with $|H|=d$, every $p$-sylowizer $S$ of $H$ in $G$ satisfies that $S\cap O^p(G)\unlhd O^p(G)$.
\end{corollary}

\begin{theorem} \label{Th7}  Let $G$ be a group, $p$ be a divisor of $|G|$ and $G_{p}\in Syl_{p}(G)$. Then $G$ is $p$-nilpotent if and only if $N_{G}(G_{p})$ is $p$-nilpotent and there exists a subgroup $P$ such that $G_{p}'\leq P\leq \Phi(G_{p})$ and $P$ has a  $p$-sylowizer $S$ satisfies that $S\cap O^p(G)$ is $S$-permutable in $G$.
\end{theorem}

\demo  The necessity is evident. Assume that $G$ is $p$-nilpotent. Then there is a normal $p$-complement $K$ such that $G=G_{p}\ltimes K$. Clearly, $O^p(G)=K$. We can pick $P=\Phi(G_{p})$. Then $P$ has a unique $p$-sylowizer $PK$ in $G$, certainly $PK\cap K\unlhd G$. Then we prove the sufficiency. If $O_{p'}(G)\neq1$, it is easy to see that the hypothesis are inherited by $G/O_{p'}(G)$. Hence $G/O_{p'}(G)$ is $p$-nilpotent by induction and so $G$ is $p$-nilpotent. Thus we may assume that $O_{p'}(G)=1$. Clearly, $G_{p}\cap S\cap O^p(G)=P\cap O^p(G)\leq P\leq \Phi(G_{p})$. By Tate's theorem (see \cite[Satz IV. 4.7]{HB}), $S\cap O^p(G)$ is $p$-nilpotent. Since $O_{p'}(G)=1$ and $S\cap O^p(G)\lhd \lhd G$, $S\cap O^p(G)=P\cap O^p(G)$. By Lemma \ref{L7}, we have $O^p(G)\leq N_{G}(P\cap O^p(G))$. Note that $G_{p}'\leq P\leq \Phi(G_{p})$, thus $G=G_{p}O^p(G)\leq N_{G}(P\cap O^p(G))$. Then by \cite[Theorem 1.5]{Y}, $G$ is $p$-nilpotent.  \qed

\begin{theorem} \label{Th8}  Let $G$ be a group and $\mathfrak{Z}=\{G_{p}, G_{p_1}, G_{p_2},\cdots, G_{p_t} \}$ be a complete set of Sylow subgroups of $G$, where $p, p_{1},\cdots, p_{t}\in \pi(G)$. Then $G$ is $p$-nilpotent if and only if $G_{p}$ has a chain of subgroups $1=P_{0}< P_{1}< \cdots < P_{n-1}< P_{n}=G_{p}$ such that $|P_i:P_{i-1}|=p$ and every $p$-sylowizer of $P_{i}$ in $G$ is $\mathfrak{Z}$-permutable in $G$ for $i=1, \cdots, n$.
\end{theorem}

\demo The necessity is evident. Since $G$ is $p$-nilpotent, there is a normal $p$-complement $K$ such that $G=G_{p}\ltimes K$. Clearly, we can pick a chain of subgroups  $1=P_{0}< P_{1}< \cdots < P_{n-1}< P_{n}=G_{p}$ of $G_{p}$ such that $|P_i:P_{i-1}|=p$ and $P_{i}\unlhd G_{p}$ for $i=1, \cdots, n$. Then $P_{i}K\unlhd G$ and $P_{i}K$ is the unique $p$-sylowizer of $P_{i}$ in $G$ by Lemma \ref{L1}, for $i=1, \cdots, n$. Next, we prove the sufficiency. Let $G$ be a counterexample with minimal order.

(1) $O_{p'}(G)=1$.

Assume that $O_{p'}(G)\neq 1$. Then $G_{p}O_{p'}(G)/O_{p'}(G)$ is a Sylow $p$-subgroup of $G/O_{p'}(G)$ and $1< P_{1}O_{p'}(G)/O_{p'}(G)< \cdots < P_{n-1}O_{p'}(G)/O_{p'}(G)< P_{n}O_{p'}(G)/O_{p'}(G)=G_{p}O_{p'}(G)/O_{p'}(G)$ is a chain of $G_{p}O_{p'}(G)/O_{p'}(G)$ such that $|P_iO_{p'}(G)/O_{p'}(G):P_{i-1}O_{p'}(G)/O_{p'}(G)|=p$ for $i=1, \cdots, n$. By the hypothesis and Lemma \ref{L5}, every $p$-sylowizer of $P_{i}O_{p'}(G)/O_{p'}(G)$ in $G/O_{p'}(G)$ a is $\mathfrak{Z}N/N$-permutable in $G/O_{p'}(G)$ for $i=1, \cdots, n$. This shows that $G/O_{p'}(G)$ satisfies the hypothesis of the theorem. Hence $G/O_{p'}(G)$ is $p$-nilpotent by the minimality of $G$, which implies that $G$ is $p$-nilpotent, a contradiction.

(2) $S_{n-1}=P_{n-1}\unlhd G$.

Let $S_{i}$ be a $p$-sylowizer of $P_{i}$ in $G$. Since $S_{i}$ is $\mathfrak{Z}$-permutable in $G$, we have $S_{i}G_{p_j}=G_{p_j}S_{i}$ for any $j=1, \cdots, t$ and $i=1, \cdots, n$. Clearly, $P_{i}$ is a Sylow $p$-subgroup of $S_{i}G_{p_j}$. Hence $G_{p_j}\leq S_{i}$ by the maximality of $S_{i}$. In particular, $G_{p_j}\leq S_{n-1}$ for any $j=1, \cdots, t$. It implies that $|G:S_{n-1}|=p^\alpha$ and so $G=G_{p}S_{n-1}$. Let $g=xy$ be an element of $G$, where $x\in S_{n-1}$ and $y\in G_{p}$. Since $P_{n-1}\unlhd G_{p}$, we have $P_{n-1}=(P_{n-1})^y\leq (S_{n-1})^y=((S_{n-1})^{x})^{y}=(S_{n-1})^g$. Hence $(S_{n-1})^g$ is a $p$-sylowizer of $P_{n-1}$ in $G$. Clearly, $G_{p_j}\leq (S_{n-1})^g$ for $j=1, \cdots, t$. It implies that $G_{p_j}\leq \cap_{g\in G}(S_{n-1})^g=(S_{n-1})_G$ for $j=1, \cdots, t$, and so $(G_{p_j})^G\leq S_{n-1}$. It follows that $O^p(G)\leq S_{n-1}$. Thus $S_{n-1}=P_{n-1}O^p(G)$. Since $P_{n-1}\unlhd G_{p}$, $S_{n-1}\unlhd G$ and $S_{n-1}$ is the unique $p$-sylowizer of $P_{n-1}$ in $G$.

Clearly, $S_{i}\leq S_{n-1}$ for any $p$-sylowizer $S_{i}$ of $P_{i}$ in $G$ and $i=1, \cdots, n-2$. Actually, if $S_{i}\nleq S_{n-1}$ for some $p$-sylowizer $S_{i}$, then there exists $p'$-element $x$ in $ S_{i}$ such that $x\notin S_{n-1}$. It implies that $P_{n-1}$ is a Sylow $p$-subgroup of $S_{n-1}\langle x \rangle$, which contradicts the maximality of $S_{n-1}$. Now, assume that $T_{i}$ is any $p$-sylowizer of $P_{i}$ in $S_{n-1}$ for $i=1, \cdots, n-2$, then by Lemma \ref{L6}, there is a $p$-sylowizer $S_{i}$ of $P_{i}$ in $G$ such that $T_{i}=S_{i}\cap S_{n-1}$ and so $T_{i}=S_{i}$ is $\mathfrak{Z}$-permutable in $G$ for $i=1, \cdots, n-2$. By Lemma \ref{L4}(2), $T_{i}$ is $\mathfrak{Z}\cap S_{n-1}$-permutable. This shows that $S_{n-1}$ satisfies the hypothesis of the theorem. Thus $S_{n-1}$ is $p$-nilpotent by the minimality of $G$. Since $O_{p'}(G)=1$, $P_{n-1}=S_{n-1}\unlhd G$.

(3) Final contradiction.

Consider the group $G/P_{n-1}$. Clearly, $G_{p}/P_{n-1}$ is a Sylow $p$-subgroup of $G/P_{n-1}$. By the hypothesis and Lemma \ref{L5}, every $p$-sylowizer of $G_{p}/P_{n-1}$ is $\mathfrak{Z}P_{n-1}/P_{n-1}$-permutable in $G/P_{n-1}$. This shows that $G/P_{n-1}$ satisfies the hypothesis of the theorem and so $G/P_{n-1}$ is $p$-nilpotent by the minimality of $G$. Let $U/P_{n-1}$ be a normal $p$-complement of $G/P_{n-1}$. Then $P_{n-1}$ is a Sylow $p$-subgroup of $U$. By the maximality of $S_{n-1}$, $U\leq S_{n-1}$ is a $p$-subgroup, a contradiction. This contradiction completes the proof. \qed

\begin{corollary} \label{Co9}
Let $G$ be a group, $p$ be a divisor of $|G|$ and $G_{p}$ a Sylow $p$-subgroup of $G$. Then $G$ is $p$-nilpotent if and only if $G_{p}$ has a chain of subgroups $1=P_{0}< P_{1}< \cdots < P_{n-1}< P_{n}=G_{p}$ such that $|P_i:P_{i-1}|=p$ and every $p$-sylowizer of $P_{i}$ in $G$ is $S$-permutable in $G$ for $i=1, \cdots, n$.
\end{corollary}

\begin{corollary} \label{Co10}
Let $G$ be a group. Then $G$ is nilpotent if and only if for every prime divisor $p$ of $|G|$ and a Sylow $p$-subgroup $G_{p}$ of $G$, $G_{p}$ has a chain of subgroups $1=P_{0}< P_{1}< \cdots < P_{n-1}< P_{n}=G_{p}$ such that $|P_i:P_{i-1}|=p$ and every $p$-sylowizer of $P_{i}$ in $G$ is $S$-permutable in $G$ for $i=1, \cdots, n$.
\end{corollary}

\begin{corollary} \label{Co11}
Let $G$ be a group. Then $G$ is $p$-nilpotent if and only if $G$ has a normal subgroup $N$  such that $G/N$ is $p$-nilpotent, and $N_{p}$ has a chain of subgroups $1=P_{0}< P_{1}< \cdots < P_{n-1}< P_{n}=N_{p}$ such that $|P_i:P_{i-1}|=p$ and every $p$-sylowizer of $P_{i}$ in $G$ is $S$-permutable in $G$ for $i=1, \cdots, n$, where $p$ is a divisor of $|N|$ and $N_{p}\in Syl_{p}(N)$.
\end{corollary}

\demo  The necessity is evident. We only need to prove the sufficiency. By the hypothesis, $G/N$ is $p$-nilpotent. Let $T/N$ be a normal $p$-complement of $G/N$. Then $N_{p}$ is a Sylow $p$-subgroup of $T$. Let $T_{i}$ be any $p$-sylowizer of $P_{i}$ in $T$. Then by Lemma \ref{L6}, $T_{i}=S_{i}\cap T$ for some $p$-sylowizer $S_{i}$ of $P_{i}$ in $G$. By the hypothesis and Lemma \ref{L3}, $T_{i}$ is $S$-permutable in $T$ for $i=1, \cdots, n$. By corollary \ref{Co9}, $T$ is $p$-nilpotent. Let $T_{p'}$ be a normal $p$-complement of $T$, then $T_{p'}\unlhd G$ and so $G$ is $p$-nilpotent.  \qed

\end{document}